\documentclass{llncs}
\usepackage{amsmath}
\usepackage{amssymb}
\begin{document}
\title{Generating matrices of highest order over a finite field}
\author{Piyasi Choudhury}
\institute{}
\date{}
\maketitle

\newtheorem{heuristic}{\bf Heuristic}
\newtheorem{algorithm}{\bf Algorithm}

\begin{abstract}
Shift registers/Primitive polynomials find applications in various branches of Mathematics, Coding Theory and Cryptography. Matrix analogues of primitive polynomials do exist. In this paper, an algorithmic approach to generating all such matrices over GF($2$) has been presented. A technique for counting all such $n \times n$ matrices over GF($2$) is also presented. The technique may be easily extended to other finite fields.
\end{abstract}


\section{Introduction}
\label{intro}

Shift registers have been widely used in the later half of the last century to generate sequences of $0$'s and $1$'s having certain other properties. A few of them are mentioned over here in brief. For a detailed account, one may refer to \cite{G67}.
\begin{enumerate}

\item{Secure and limited-access code generators}
\begin{enumerate}
\item{Encipherment}
\item{Privacy encoding}
\item{Multiple address coding}
\end{enumerate}

\item{Efficiency code generators}
\begin{enumerate}
\item{Error-correcting codes}
\item{Signals recoverable through noise}
\end{enumerate}

\item{Prescribed property generators}
\begin{enumerate}
\item{Prescribed period generators}
\item{Prescribed sequence generators}
\end{enumerate}

\item{Mathematical models}
\begin{enumerate}
\item{Random bit generators}
\item{Finite state machines}
\item{Markov processes}
\end{enumerate}

\end{enumerate}

In \cite{G67}, it has been shown that each state of an $n$-tube shift register can be thought of as an $n$-dimensional vector. The shift register acts as a linear operator and changes each state into the next. This linear operator operating on $n$-dimensional vectors is represented by a $n \times n$ matrix. The characteristic equation of this matrix is essentially the same as the characteristic polynomial of the shift register. In this paper, only matrices are dealt with.

The matrices under consideration belong to the last two application areas.

All the $n \times n$ matrices of highest order over GF($2$) have been counted with detailed description of each step. Begin with the property of the characteristic polynomial of these matrices and note that these polynomials are primitive polynomials. Also note that these matrices are the primitive elements of the general linear group $GL_n(GF(2))$,  conjugate to each other and belong to an equivalence class where the the equivalence relation is conjugacy. The cardinality of this equivalence class is calculated from elementary abstract algebra. Once all the above steps are complete, it remains an easy task to calculate the  total number of such matrices. A method for generating all such matrices is presented next.

\section{Preliminaries: A few definitions and results from Abstract Algebra}
At the outset, a few results are presented. For the proof and other details, one may refer to \cite{G67,LN86,MS77}
\begin{theorem}
The maximum order of an $n \times n$ matrix $A$ over $GF(2)$ is $2^n - 1$.
\end{theorem}

\begin{theorem}
An $n \times n$ matrix $A$ over $GF(2)$ has maximum order $2^n - 1$ iff its characteristic polynomial is a primitive polynomial of degree $n$.
\end{theorem}

\begin{example}
The matrix 
\begin{math}\bordermatrix{
&      &      &      \cr
& 0    & 1    & 0    \cr          
& 0    & 0    & 1    \cr          
& 1    & 1    & 0 \cr}
\end{math}

has the characteristic polynomial $f(x) = x^3 + x + 1$, which is a primitive polynomial of degree $3$ over $GF(2)$. Hence its order is $2^3 - 1 = 7$
\end{example}

\begin{theorem}
There exist equal number of matrices corresponding to each of the primitive polynomials of degree $n$.
\end{theorem}

\begin{theorem}
There are $\frac {\Phi(2^n - 1)}{n}$ number of primitive polynomials of degree $n$ where $\Phi(n)$ is Euler's Totient function.
\end{theorem}

In order to count the $n \times n$ matrices of order $2^n - 1$, it will suffice to count the matrices whose chararcteristic polynomial is a given primitive polynomial, say, $f(x)$.

Since the characteristic polynomial of A is a primitive polynomial, the minimal polynomial of A is also the same primitive polynomial, $f(x)$.

\textbf{Definition: Companion Matrix}\\
In general, the \textit{companion matrix} of a monic polynomial $f(x) = a_0 + a_1x+ \ldots +a_{n-1}x^{n-1}+x^n$ of positive degree $n$ over a field is defined to be the $n \times n$ matrix

\begin{math}\bordermatrix{
&      &      &			&		 &   &  \cr
& 0    & 0    & 0   &\ldots	 & 0    & -a_0 \cr
& 1    & 0    & 0   &\ldots		 & 0    & -a_1 \cr
& 0    & 1    & 0   &\ldots		 & 0    & -a_2 \cr
&\vdots    &\vdots  &\vdots   &\vdots	 &\vdots    &\vdots \cr
& 0    & 0    & 0   &\ldots		 & 1   & -a_{n-1}    \cr}
\end{math}

\begin{theorem}
All the matrices whose characteristic polynomial and the minimal polynomial are same and equal to $f(x)$, (where $f(x)$ is a primitive polynomial), are similar or conjugate to each other.
\end{theorem}

The proof of the above theorem is based on the fact that all such matrices are similar to the companion matrix of the polynomial $f(x)$. A direct application of the \textit{Rational Form Theorem} proves the result easily.

\begin{example}
For $n = 3$, there are two such polynomials. In other words, there are two primitive polynomials of degree $3$ as follows:
$$f_1(x) = x^3 + x + 1$$ and $$f_2(x) = x^3 + x^2 + 1$$

The companion matrix of $f_1(x)$ is 
\begin{math}\bordermatrix{
&      &      &      \cr
& 0    & 0    & 1    \cr          
& 1    & 0    & 1    \cr          
& 0    & 1    & 0 \cr}
\end{math}

The companion matrix of $f_2(x)$ is 
\begin{math}\bordermatrix{
&      &      &      \cr
& 0    & 0    & 1    \cr          
& 1    & 0    & 0    \cr          
& 0    & 1    & 1    \cr}
\end{math}\\\\
The above theorem asserts that all the matrices having the characteristic (here, it is the same as the minimal polynomial) $f_1(x)$ are similar. Also all the matrices corresponding to $f_2(x)$ are imilar.

\end{example}

It may be noted that the matrices we are dealing with are elements of the general linear group $GL_n(GF(2))$, the group of non-singular $n \times n$ matrices over $GF(2)$.

\textbf{Definition: General Linear Group over $GF(2)$ and its order}\\
The set of all $n \times n$ invertible matrices over a field $F$ forms a group with respect to matrix multiplication.This group is called the general linear group of degree $n$ over $F$ and is denoted by $GL_n(F)$. The identity element of $GL_n(F)$ is the identity matrix and the inverse element of $A$ is $A^{-1}$. It may be noted that $GL_n(F)$ is not an abelian group for $n \geq 2$. In the following treatment, we shall choose $F$ as $GF(2)$ and we shall refer to $GL_n(GF(2))$ as $G$ henceforth.

The number of elements in $G$ is given by 
$$ \prod^{n - 1}_{i = 0}(2^n - 2^i)$$

The proof of the above is easy to see:
\begin{enumerate}
\item The first column of a non-singular matrix must not be the $0$ vector,thus there are $2^n-1$ many to form the first column.
\item The $i$-th column must not be a linear combination of the previous $i-1$-columns, thus there are only $(2^n - 2^{i-1})$ many choices.
\end{enumerate}

So $$|G| = \prod^{n}_{i = 1}(2^n - 2^{i-1}) = \prod^{n - 1}_{i = 0}(2^n - 2^i)$$

Let $G$ be a group.\\ 
\textbf{Definition: Conjugacy}\\
If $a,b \in G$, then b is said to be a conjugate of $a$ in $G$ if there exists an element $c \in G$ such that $b = c^{-1}ac$

Elements $a,b \in G$ related as above are called conjugate.

\begin{theorem}
Conjugacy is an equivalence relation on $G$
\end{theorem}

One can partition $G$ into disjoint equivalence classes $Cl(a)$.\\

\textbf{Definition: Centraliser}\\
If $a \in G,~ define ~N(a) = \{x \in G| xa = ax\}$. It is easy to verify that $N(a)$ is a subgroup of $G$. $N(a)$ is usually called the centraliser of $a$ in $G$.\\

\begin{theorem}
If $a \in G$, $|N(a)| = [g:N(a)] = \frac{|G|}{|N(a)|}$
\end{theorem}

Consider $A$ as a linear operator on the finite-dimensional vector space $V$ over $GF(2)$.\\
\textbf{Definition: Cyclic vector}\\
A vector $\alpha \in V$ is called a cyclic vector for $A$ if the vectors $A^k \alpha, k = 0,1,\ldots,$ span $V$.

\begin{theorem}
Let $A$ be a linear operator on the finite-dimensional vector space $V$. Then $A$ has a cyclic vector if and only if the characteristic and minimal polynomials for $A$ are identical.
\end{theorem}

\textbf{Definition: Cyclic matrix}\\
If the characteristic polynomial of a matrix is the product of distinct irreducible factors, the matrix is said to be a cyclic matrix.

Clearly, the characteristic polynomial and the minimal polynomial are the same for cyclic matrices.

\section{The counting: step by step}
\begin{theorem}
If the characteristic polynomial of a matrix $f(x)$ be irreducible, it may be shown that any non-zero polynomial in $A$ is non-singular.
\end{theorem}

(The above theorem may be used since we are dealing with primitive polynomials of degree $n$ that are necessarily irreducible)

\begin{proof}
Let $g(A)$ be such a polynomial. Since $f(x)$ is irrducible, $g(x)$ and $f(x)$ are relatively prime.

From Euclid's algorithm,$$1 = p(x)g(x)+ q(x)f(x)$$ for some polynomials $p(x)$ and $q(x)$.

Replacing $x$ by $A$ and $1$ by $I$, the following matrix identity is obtained:
$$I = p(A)g(A)+ q(A)f(A) = p(A)g(A)$$ since $f(A) = 0$, by definition of the minimal polynomial.

Hence $g(A)$ is invertible and therefore non-singular.
\end{proof}

\begin{theorem}
If $A$ is a cyclic matrix and $B$ commutes with $A$, then $B$ is a polynomial in $A$
\end{theorem}

\begin{proof}
Since $A$ is a cyclic matrix, there exists a cyclic vector $\alpha$ such that the set $S = {v,Av, \ldots,A^{n-1}v}$ is a basis of $V$.

Consider a matrix $B$ that commutes with $A$. The vector $B\alpha$ can be expressed as a linear combination of the vectors in $S$ since $S$ is a basis.

$B\alpha = a_0 + a_1A\alpha + \ldots + a_{n-1}A^{n-1}\alpha = g(A)\alpha$

where $g(A) = a_0 + a_1A + \ldots + a_{n-1}A^{n-1}$ is a polynomial in $A$ (of degree $(n-1)$).

As $B$ commutes with $A$, it commutes with any power of $A$ and also with $g(A)$.

Obviously, $A$ also commutes with $g(A)$.

If the linear operator $B$ is applied to any vector $(A^k)\alpha$ in $S$, one gets
$B(A^k)\alpha = (A^k)B\alpha = (A^k)g(A)\alpha = g(A)(A^k)\alpha$

Finally consider a matrix $C = [\alpha | A\alpha | \ldots | A^{n-1}\alpha]$

The above result shows that the $k$-th column of $BC$ is equal to the $k$-th column of $g(A)C$ for all $k = 1,2, \ldots , n$

which means $BC = g(A).C$.

Since $C$ is non-singular, so $B = g(A)$.
\end{proof}

\textbf {\large{The final result:}}\\
\begin{theorem}
The number of non-zero polynomials of degree less than $n$ is $2^n - 1$
\end{theorem}

\begin{proof}
Consider the $(n-1)$-th degree polynomial $c_0 + c_1x + \ldots + c_{n-1}A^{n-1}$. There are $n$ coefficients and each of them can assume only two values, zero or one.

Hence there are totally $2^n -1$ non-zero polynomials of degree less than $n$.
\end{proof}

Since $B$ is a polynomial of degree less than $n$ in $A$, $|N(A)| = 2^n - 1$

So the number of matrices (say,$A$) with characteristic polynomial $f(x)$ is equal to $$|Cl(A)| = [G:N(A)] = \frac{|G|}{|N(A)|} = \frac {\prod^{n-1}_{i=0}(2^n - 2^i)}{2^n-1} = \prod^{n-1}_{i=1}(2^n - 2^i)$$

The above is the number of matrices per primitive polynomial.

Since there are $\frac {\Phi(2^n - 1)}{n}$ number of primitive polynomials, so there are totally $$\prod^{n-1}_{i=1}(2^n - 2^i)\frac {\Phi(2^n - 1)}{n}$$ number of matrices of order $2^n -1$.\\\\

\section{The Method}
\textbf{Given:} A primitive polynomial of degree $n$. Let us denote it by $f(x)$.\\
\textbf{Steps:}\\
\begin{enumerate}
\item Calculate the companion matrix of $f(x)$. Call it $A$.

\item Calculate the centraliser $H$ of $A$ in $GL_n(GF(2))$, which is nothing but the collection of all the polynomials in $A$ of degree less than $n$.

\item Find the coset decomposition of $GL_n(GF(2))$ with respect to $H$. Let $\frac{|G|}{|H|} = k$. Let the cosets be $c_0, c_1, \ldots, c_{k-1}$ and $|H| = |c_0| =  |c_1| =  \ldots =  |c_{k-1}| = m$(say).

\item Choose one matrix each from these cosets. Denote the matrix chosen from the $i$-th coset $c_i$ as $M_i$ for $i = 0,1,\ldots,k-1$.

\item Calculate the $k$ number of conjugates of $A$ using $M_i$, i.e. $M_iAM_i^{-1}$ for $i = 0,1,\ldots,k-1$.

\end{enumerate}

It is clear that two matrices belonging to the same right coset of $H$ in $G$ will yield the same conjugate of $A$ and two matrices belonging to different right cosets of $H$ in $G$ will yield different conjugates of $A$.

Thus, if the coset decomposition of $G$ is known, calculating the matrices similar (conjugate) to $A$ is easy. All these matrices are of maximum order.\\\\

\textbf{\Large{Illustration for $n = 3$}}
\begin{enumerate}
\item \textbf{\large{Notation:}} The following matrix 

\begin{math}\bordermatrix{
&      &   &  \cr
& a_{00}    & a_{01}    & a_{02} \cr
& a_{10}    & a_{11}    & a_{12} \cr
& a_{20}    & a_{21}    & a_{22}  \cr}
\end{math}

is represented by the integer $a_{00}2^0 + a_{01}2^1 + a_{02}2^2 + a_{10}2^3 +  a_{11}2^4 + a_{12}2^5 + a_{20}2^6 + a_{21}2^7 + a_{22}2^8$. It may be noted that there are $512$ number of $3 \times 3$ matrices and they will be represented by the integers $0,1,\ldots,511$.

\item \textbf{\large{Coset decomposition of $GL_3(GF(2))$}}

Note that $|GL_3(GF(2))| = 168$.

Consider the matrix $A = $
\begin{math}\bordermatrix{
&      &      &      \cr
& 0    & 0    & 1    \cr          
& 1    & 0    & 1    \cr          
& 0    & 1    & 0 \cr}
\end{math}

This is the companion matrix of the polynomial $x^3+x+1$.

According to our notation, $A$ is represented as $172$.
Construct the subgroup $H$ by the centraliser of $A$, which consists of all the polynomials in $A$ of degree less than $3$. So $|H| = 2^3 - 1 = 7$. $H = \{106,157,247,273,379,396,486\}$. The number of cosets of $H$ in $GL_3(GF(2))$ is $\frac{|GL_3(GF(2))|}{|H|} = \frac{168}{7} = 24$. Now $H$ itself is one of the cosets. So there are $23$ other cosets. Let us denote them by $c_0,c_1, \ldots, c_{22}$.\\

$c_0 = \{84,169,253,346,270,499,423\}$\\
$c_1 = \{85,161,244,282,335,443,494\}$\\
$c_2 = \{86,185,239,474,396,355,309\}$\\
$c_3 = \{87,177,230,410,461,299,380\}$\\
$c_4 = \{92,233,181,339,271,442,486\}$\\
$c_5 = \{93,225,188,275,334,498,431\}$\\
$c_6 = \{94,249,167,467,397,298,372\}$\\
$c_7 = \{95,241,174,403,460,354,317\}$\\
$c_8 = \{98,281,379,236,142,501,407\}$\\
$c_9 = \{102,313,351,492,394,213,179\}$\\
$c_{10} = \{103,305,342,428,459,157,250\}$\\
$c_{11} = \{106,345,307,229,143,444,470\}$\\
$c_{12} = \{107,337,314,165,206,500,415\}$\\
$c_{13} = \{110,377,279,485,395,156,242\}$\\
$c_{14} = \{111,369,286,421,458,212,187\}$\\
$c_{15} = \{114,409,491,254,140,359,277\}$\\
$c_{16} = \{115,401,482,190,205,303,348\}$\\
$c_{17} = \{116,425,477,382,266,215,163\}$\\
$c_{18} = \{117,417,468,318,331,159,234\}$\\
$c_{19} = \{122,473,419,247,141,302,340\}$\\
$c_{20} = \{123,465,426,183,204,358,285\}$\\
$c_{21} = \{124,489,405,375,267,158,226\}$\\
$c_{22} = \{125,481,412,311,330,214,171\}$\\

Now choose one coset each from the above $23$ cosets and calculate the conjugate of $A$ using that matrix. Since all the matrices belonging to the same coset will yield the same conjugate, choosing any one of the matrices from each coset will give the same result. Finally we shall get a list of 24 matrices which are similar (conjugate) to the matrix $A$. The list is as follows:\\
$\{95,335,187,485,442,500,102,142,172,226,106,204,115,397,157,355,370,\\412,247,431,253,491,382,478\}$.\\

All these $24$ matrices have the same characteristic/minimal polynomial, which is $x^3+x+1$, and all of them are of maximum order, that is $7$.

Next consider the matrix $B = $
\begin{math}\bordermatrix{
&      &      &      \cr
& 0    & 0    & 1    \cr          
& 1    & 0    & 0    \cr          
& 0    & 1    & 1    \cr}
\end{math}

This is the companion matrix of the polynomial $x^3+x^2+1$.

According to our notation, $B$ is represented as $396$. Following similar lines as above, we construct the subgroup consisting of the polynomials of degree less than three and then the $23$ cosets. Here,\\
$H = \{106,157,247,273,379,396,486\}$ \\ and \\
$c_0 = \{84,225,181,459,415,298,382\}$\\ 
$c_1 = \{85,233,188,395,478,354,311 \}$\\
$c_2 = \{86,241,167,331,285,442,492 \}$\\
$c_3 = \{87,249,174,267,348,498,421 \}$\\
$c_4 = \{92,161,253,394,470,299,375 \}$\\
$c_5 = \{93,169,244,458,407,355,318 \}$\\
$c_6 = \{94,177,239,266,340,443,485 \}$\\
$c_7 = \{95,185,230,330,277,499,428 \}$\\
$c_8 = \{98,337,307,461,431,156,254 \}$\\
$c_9 = \{99,345,314,397,494,212,183 \}$\\
$c_{10} = \{102,369,279,205,171,444,474\}$\\
$c_{11} = \{103,377,286,141,234,500,403\}$\\
$c_{12} = \{107,281,370,460,423,213,190\}$\\
$c_{13} = \{110,305,351,140,226,445,467\}$\\
$c_{14} = \{111,313,342,204,163,501,410\}$\\
$c_{15} = \{114,465,419,335,317,158,236\}$\\
$c_{16} = \{115,473,426,271,380,214,165\}$\\
$c_{17} = \{116,481,405,207,187,302,346\}$\\
$c_{18} = \{117,489,412,143,250,358,275\}$\\
$c_{19} = \{122,401,491,270,372,159,229\}$\\
$c_{20} = \{123,409,482,334,309,215,172\}$\\
$c_{21} = \{124,417,477,142,242,303,339\}$\\
$c_{22} = \{125,425,468,206,179,359,282\}$\\

Choosing one matrix from each coset and then calculating the conjugate using that matrix finally yields a list of 24 matrices which are similar (conjugate) to the matrix $B$. The list is as follows:\\

$\{244,426,229,171,334,94,156,354,99,141,396,114,492,250,486,190,207,\\111,379,477,415,375, 499,445\}$

All these $24$ matrices have the same characteristic/minimal polynomial, which is $x^3+x+1$, and all of them are of maximum order, that is $7$.

\item \textbf{\large{The matrices corresponding to the polynomial $x^3+x+1$}}
$\{95,335,187,485,442,500,102,142,172,226,106,204,115,397,157,355,370,\\412,247,431,253,491,382,478\}$.\\

\item \textbf{\large{The matrices corresponding to the polynomial $x^3+x^2+1$}}
$\{244,426,229,171,334,94,156,354,99,141,396,114,492,250,486,190,207,\\111,379,477,415,375, 499,445\}$

\end{enumerate}

\section{Conclusion}
It is clear from the above example for $n = 3$ and $GF(2)$ that one can find all the $n \times n$ matrices of maximum order, i.e., $|2^n-1|$, given the coset decomposition of the General Linear Group of order $n$. These matrices are equivalent to shift registers, as pointed out earlier. If one knows the entire pool of matrices that exhibit the \textit{primitive} behaviour, one can choose any one from the pool and use it for an application where \textit{random} behaviour is sought for. If the matrix is thought of as a linear operator/transformation (operating on a non-zero vector of length $n$ called the \textit{state}), it is immediate that the resulting sequence of states will comprise of all possible $2^n-1$ states and also this is the cycle length of the \textit{state diagram}. However, the order of the states in each case will be different (if one starts with a new matrix each time), thus providing some notion of randomness.

\end{document}